\documentclass{amsart}

\usepackage{amssymb,graphicx, color}
\usepackage{xypic}
\usepackage[all]{xy}
\usepackage{lscape}
\usepackage{longtable}

\newtheorem{thm}{Theorem}[section]

\newtheorem{lem}[thm]{Lemma}

\theoremstyle{definition}
\newtheorem{defn}[thm]{Definition}

\newtheorem{rem}[thm]{Remark}

\newtheorem{defn-thm}[thm]{Definition--Theorem}  
\newtheorem{defn-lem}[thm]{Definition--Lemma}  

\theoremstyle{remark}

\renewcommand{\c}[0]{{\mathbb C}}

\newcommand{\p}[0]{{\mathbb P}}

\newcommand{\rank}[0]{\operatorname{rank}}

\newcommand{\C}{\mathbb{C}}

\def\loccoh#1.#2.#3.#4.{H^{#1}_{#2}(#3,#4)}

\DeclareMathAlphabet{\mathchanc}{OT1}{pzc}%
                                {m}{it}

\numberwithin{equation}{section}

\begin{document}
\bibliographystyle{amsalpha}

\title[Waring decompositions and identifiability via Bertini and Macaulay2...]{Waring decompositions and identifiability \\ via Bertini and Macaulay2 software}

\author[Elena Angelini]{Elena Angelini}
\address[Elena Angelini]{Dipartimento di Ingegneria dell'Informazione e Scienze Matematiche, Universit\`a di Siena, Via Roma 56, 53100 Siena, Italia}
\email{elena.angelini@unisi.it}

\begin{abstract}
Starting from our previous papers \cite{AGMO} and \cite{ABC}, we prove the existence of a non-empty Euclidean open subset whose elements are polynomial vectors with $4$ components, in $ 3 $ variables, degrees, respectively, $ 2,3,3,3 $ and rank $ 6 $, which are not identifiable over $ \c $ but are identifiable over $ \mathbb{R} $. This result has been obtained via computer-aided procedures suitably adapted to investigate the number of Waring decompositions for general polynomial vectors over the fields of complex and real numbers. Furthermore, by means of the Hessian criterion (\cite{COV}), we prove identifiability over $ \C $ for polynomial vectors in many cases of sub-generic rank.


\noindent \emph{Mathematics~Subject~Classification~(2010):}
15A69, 14P05, 14Q99, 13P05, 65H10.
\end{abstract}

\maketitle

\section{Introduction}\label{sec:intr}

Waring identifiability is a quite rare phenomenon, investigated between the XIX\textsuperscript{th} and the XX\textsuperscript{th} century. 

Precisely, given a general homogeneous polynomial $ p \in \C[x_{0}, \ldots, x_{n}]_{d} $, we say that $ p $ is \emph{Waring identifiable} if it admits a unique decomposition (up to reordering and rescaling) as linear combination of $ d $\textsuperscript{th} powers of elements in $ \C[x_{0}, \ldots, x_{n}]_{1} $. When this happens, the decomposition provides
a canonical form of $p$. An expression as described above is called a \emph{Waring decomposition} of $ p $, in honor of Waring statement (1770) in Number Theory. The minimum number of summands $ k $ in a Waring decomposition is the \emph{Waring rank} of $ p $. 

From the mathematical literature, we know that a general $ p $ has a unique Waring decomposition with rank $ k $ if $(n,d,k) \in \{(1,2t+1,t+1),(2,5,7),(3,3,5)\}$, $ t \in \mathbb{N} $. The first and third cases are due to Sylvester (\cite{Sy}) while the second to Hilbert (\cite{Hi}). The third one is the so called \emph{Sylvester's Pentahedral Theorem}. Recently, in \cite{GM}, it has been proved that these are the only examples. 

This phenomenon is very important in applications, e.g. in Blind Source Separation, Phylogenetic and Algebraic Statistics, we refer to \cite{Lan} for these aspects. It has been developed also in the more general setting of tensors. \\
\indent When Waring identifiability fails, we can attempt to recover it by introducing the simultaneous decomposition of more forms with the same number of summands or by requiring real identifiability. We describe these two new points of view and how they complement one another.

Concerning the first aspect, we say that a \emph{simultaneous Waring decomposition} of a \emph{polynomial vector} $ f= (f_1, \ldots, f_r) $ of general forms $ f_{j} \in \C[x_{0}, \ldots, x_{n}]_{d_{j}} $, with $ d_{1} \leq \ldots \leq d_{r} $, is a way of writing its components $ f_1, \ldots, f_r $ as a linear combination of powers of the same linear forms. The notions of Waring identifiability and rank extends naturally to the case $ r>1 $. 

For instance, if $ r= d_{1} = d_{2} = 2 $, then we deal with the simultaneous diagonalization (by congruence) of a pair of quadrics in $ \p_{\C}^n $: $ f = (f_{1},f_{2}) $ admits a unique Waring decomposition with rank $ k=n+1 $, the so called \emph{Weierstrass canonical form} (\cite{We}), if and only if the discriminant of the pencil $ \left<f_{1},f_{2}\right> $ does not vanish. In particular, if $ f_{1}$ is the Euclidean quadratic form, then we get the diagonalization of $ f_{2} $ with orthogonal summands, with respect to $f_{1}$; over $ \mathbb{R} $ this is possible for any $ f_{2} $ thanks to the \emph{Spectral Theorem}. Other known cases of Waring identifiability are the following: 
$$ (r,n,d_1,\ldots, d_{r},k) = \left\{\begin{array}{l}
\left(r,1,d_{1}, \ldots, d_{r},\left \lceil {\frac{1}{1+r}}\sum_{j=1}^{r}{1+d_{j} \choose d_{j}} \right \rceil\right), \,d_{1}+1 \geq k,\, \textrm{\cite{CR}}\\
(4,2,2,2,2,2,4), \textrm{Veronese \cite{Ve}}\\
(2,2,2,3,4), \textrm{Roberts \cite{Ro}} \\
(3,2,3,3,4,7), \textrm{\cite{AGMO}}
\end{array}\right. $$

Regarding the second aspect, we can ask if a general form or polynomial vector, with real coefficients, which is not Waring identifiable, admits a unique decomposition with real linear forms and real scalars. When this happens, we speak about \emph{real identifiability}. In this direction, in \cite{ABC} we proved that, when $ r=1 $, real identifiability holds in non-trivial Euclidean open subsets in the space of homogeneous polynomials with $(n,d,k)\in \{(2,7,12);(2,8,15)\}$. We notice that, in these examples, the number of Waring decompositions for a general form is, respectively, 5 (\cite{DS}) and 16 (\cite{RS}).\\
\indent By combining the simultaneous decomposition problem with the real point of view, we are able to prove the following:
\begin{thm}\label{thm:23330}
There exists a non-trivial Euclidean open subset in the space of real polynomial vectors with $4$ components, in $ 3 $ variables, degrees, respectively, $ 2,3,3,3 $ whose elements have rank $ 6 $ and are identifiable over $ \mathbb{R} $ but not over $ \C $.
\end{thm}
The real identifiability stated in Theorem \ref{thm:23330} arises from the computational analysis described in section \ref{sss:cap}, realized via the software for \emph{Numerical Algebraic Geometry} Bertini (\cite{Be}, \cite{BHSW}), suitably coordinated with Matlab, in the spirit of \cite{AGMO} and \cite{ABC}. Our technique is based on \emph{homotopy continuation} and \emph{monodromy loops} (\cite{HOOS}).

By means of a different computational approach developed with the software Macaulay2 (\cite{M2}), which we outline in section \ref{sss:hp}, we are also able to produce many cases of simultaneous Waring identifiability from the complex side, under the assumption of \emph{sub-generic rank}. Our method is based on the \emph{Hessian criterion}, originally introduced in \cite{COV} to study the analogous problem in the case of tensors. Our results in this direction are collected in section \ref{sss:res}.

The paper is organized as follows. In section $ 2 $ we recall main definitions and known results about Waring decomposition, rank and identifiability, focusing on the complex and real side; moreover we introduce the simultaneous Waring setting, for a detailed description of which we refer to \cite{AGMO}. Section $ 3 $ is devoted to our main results and is divided into $ 2 $ parts: in section \ref{ss:NAG} we describe the computational approach and the proof for our Theorem \ref{thm:23330}; in section \ref{sec:hes} we focus on sub-generic rank cases, presenting our procedure and the identifiable cases detected with it.

\section{Preliminaries}

\subsection{Waring decomposition and rank}\label{ss:dr}$ \quad $\\

Let $\mathbb{F}$ be either the complex field $ \mathbb{C} $ or the real field $ \mathbb{R} $. Let $ (n, d) \in \mathbb{N}^{2} $ and let $ \mathbb{F}[x_{0}, \ldots, x_{n}]_{d} $ be the space of homogeneous polynomials, called for simplicity \emph{forms}, of degree $ d $ in $ n+1 $ variables over $\mathbb{F}$.
\begin{defn}\label{def:Wardec}
A \emph{Waring decomposition} over $ \mathbb{F} $ of $ p \in \mathbb{F}[x_{0}, \ldots, x_{n}]_{d} $ is given by linear forms $ \ell_{i} \in \mathbb{F}[x_{0}, \ldots, x_{n}]_{1} $ and scalars $ \lambda_{i} \in \mathbb{F} - \{0\}  $, $ i \in\{1, \ldots, k\} $, such that
\begin{equation}\label{eq:Wardec}
p = \displaystyle{\sum_{i=1}^{k} \lambda_{i}\ell_{i}^d}. 
\end{equation}
\end{defn}
\begin{defn}\label{def:rank}
The minimal $ k $ appearing in (\ref{eq:Wardec}) is called the \emph{Waring} \emph{rank} of $ p $ over $ \mathbb{F} $. For simplicity, throughout the paper, we speak about rank.
\end{defn}
\begin{rem}
Any summand in (\ref{eq:Wardec}) has rank $ 1 $ over $ \mathbb{F} $.
\end{rem}
\begin{defn}
A \emph{typical rank} over $ \mathbb{F} $ for forms in $ \mathbb{F}[x_{0}, \ldots, x_{n}]_{d} $ is any $ k $ such that the set of forms having rank $ k $ has positive Euclidean measure.
\end{defn}
\begin{defn}
There exists a unique typical rank over $ \mathbb{C} $ for forms in $ \mathbb{C}[x_{0},\ldots,x_{n}]_{d} $, which we call the \emph{generic rank} for that space of forms. 
\end{defn}
One of the most important results for the complex case is the following:
\begin{thm}[Alexander-Hirschowitz, 1995, \cite{AH}]\label{thm:AH}
The expression for the generic rank of $ \mathbb{C}[x_{0},\ldots,x_{n}]_{d} $ is $ k_{g} = \left \lceil{{n+d \choose d}\frac{1}{n+1}} \right \rceil$, except for
\begin{itemize}
\item $ (n,2) $, where $ k_{g} = n+1 $;
\item $ (n,d) \in \{(2,4),(3,4),(4,3),(4,4)\} $, where $ k_{g} = \left \lceil{{n+d \choose d}\frac{1}{n+1}} \right \rceil + 1$. 
\end{itemize}
\end{thm}
\begin{rem}
As a consequence of Theorem \ref{thm:AH}, if $ {n+d \choose d}\frac{1}{n+1} \in \mathbb{N}$ then the general form in $ \mathbb{C}[x_{0},\ldots,x_{n}]_{d} $ of generic rank has finitely many minimal Waring decompositions, with some of the exceptions above. 
\end{rem}
According to \cite{COV}, we give the following:
\begin{defn}\label{def:subgen}
Forms in $ \mathbb{C}[x_{0},\ldots,x_{n}]_{d} $ of rank $ k < \left\lceil{n+d \choose d}\frac{1}{n+1}\right\rceil $ are of \emph{sub-generic rank}.
\end{defn}
\begin{rem}\label{rem:c&r}
It is well known that it is possible to have more than one typical rank over $ \mathbb{R} $ for forms in $ \mathbb{R}[x_{0},\ldots,x_{n}]_{d} $. The smallest typical rank over $ \mathbb{R} $ coincides with the generic rank over $ \mathbb{C} $, \cite{BT}. Any rank over $ \mathbb{R} $ between the minimal typical rank and the maximal typical rank is also typical and it is an open problem to determine an expression for the expected maximal typical rank over $ \mathbb{R} $ in the general case, since only partial results have been obtained, \cite{BBO}.
\end{rem}

\subsection{Classical Waring identifiability}

\begin{defn}\label{def:ident}
A rank-$ k $ form $ p \in \mathbb{F}[x_{0}, \ldots, x_{n}]_{d} $ is \emph{Waring identifiable} over $ \mathbb{F} $ if the presentation (\ref{eq:Wardec}) is unique up to a permutation and scaling of the summands. 
\end{defn} 

\subsubsection{Complex identifiability}\label{sss:ci}$ \quad $\\

Being Waring identifiable is rare for forms of generic rank over $ \mathbb{C} $. Namely, as shown in \cite{GM}, this happens only in the classically known cases listed in Table $ 1 $. \newpage
\begin{table}[h]
\begin{center}
\begin{tabular} {c |c |c }
{\bf Space of forms} & {\bf Generic rank}  & {\bf Ref.} \\
$ \mathbb{C}[x_{0},x_{1}]_{2t+1} $ & $t+1$ & \cite{Sy} \\
$ \mathbb{C}[x_{0},x_{1},x_{2}]_{5} $ & $7$ & \cite{Hi} \\
$ \mathbb{C}[x_{0},x_{1},x_{2},x_{3}]_{3} $ & $5$ & \cite{Sy} \\
\end{tabular}\vspace{0.2cm}\caption{Waring identifiable cases of generic rank over $ \mathbb{C} $}
\end{center}
\end{table}

On the contrary, Waring identifiability over $ \mathbb{C} $ is expected for general forms of sub-generic rank. In this direction, in \cite{COV1} it is proved that, if $ d \geq 3 $, then there are no exceptions besides the ones appearing in Table $ 2 $.
\begin{table}[h]
\begin{center}
\begin{tabular} {c |c |c}
{\bf Space of forms} & {\bf Sub-generic rank} & {\bf Ref.} \\
$ \mathbb{C}[x_{0},x_{1},x_{2}]_{6} $ & $9$  & \cite{AC},\cite{CC} \\
$ \mathbb{C}[x_{0},x_{1},x_{2},x_{3}]_{4} $ & $8$  & \cite{CC},\cite{M}\\
$ \mathbb{C}[x_{0},x_{1},x_{2},x_{3},x_{4},x_{5}]_{3} $ & $9$  & \cite{COV1},\cite{RV}\\
\end{tabular}\vspace{0.2cm}\caption{Waring unidentifiable cases of sub-generic rank over $ \mathbb{C} $}
\end{center}
\end{table}

\subsubsection{Real identifiability}\label{sss:ri}$ \quad $ \\

Recent interest has been devoted to the real case (\cite{COV2}, \cite{CLQ}, \cite{DDL}), which is very important in applications. 

Consider a general $ p \in \mathbb{R}[x_{0}, \ldots, x_{n}]_{d} $ such that the rank over $ \mathbb{C} $ equals the rank over $ \mathbb{R} $. If $ p $, seen as an element of $ \mathbb{C}[x_{0}, \ldots, x_{n}]_{d} $, is Waring identifiable over $ \mathbb{C} $, then, necessarily, $ p $ is Waring identifiable over $ \mathbb{R} $. Namely, with the above assumption, the unique complex decomposition of $ p $ is actually completely real. 

At this point a natural question arises: if $ p $ is not Waring identifiable over $ \mathbb{C} $, what can one say about the Waring identifiability of $ p $ over $ \mathbb{R} $? Recently, in \cite{ABC} we proved that real identifiability holds in non-trivial Euclidean open subsets of the spaces of forms collected in Table $ 3 $ (in the second column, \emph{g} stands for generic and \emph{sg} for sub-generic).

\begin{table}[h]
\begin{center}
\begin{tabular} {c |c |c |c }
{\bf Space of forms} & {\bf Rank}  & {\bf Dec. over $ \mathbb{C} $} & {\bf Ref.} \\
$ \mathbb{C}[x_{0},x_{1},x_{2}]_{7} $ & $12$ (g) & 5 & \cite{DS} \\
$ \mathbb{C}[x_{0},x_{1},x_{2}]_{8} $  & $15$ (g) & 16 & \cite{RS} \\
$ \mathbb{C}[x_{0},x_{1},x_{2}]_{6} $ & $9$ (sg) & 2 & \cite{AC},\cite{CC} \\
$ \mathbb{C}[x_{0},x_{1},x_{2},x_{3}]_{4} $ & $8$ (sg) & 2 & \cite{CC},\cite{M}\\
$ \mathbb{C}[x_{0},x_{1},x_{2},x_{3},x_{4},x_{5}]_{3} $ & $9$ (sg) & 2 & \cite{COV1},\cite{RV}\\
\end{tabular}\vspace{0.2cm}\caption{Unidentifiable cases over $ \mathbb{C} $, identifiable over $ \mathbb{R} $ in open sets}
\end{center}
\end{table} 

Moreover in \cite{ABC} we showed that, if there are infinitely many minimal Waring decompositions over $ \mathbb{C} $, then it is not possible to find Euclidean open sets of forms in which real identifiability holds. 

Therefore, requiring real identifiability of forms is an efficient method to recover Waring identifiability when we deal with finitely many decompositions over $ \mathbb{C} $.

\subsection{Simultaneous Waring setting}\label{sss:sci} $ \quad $\\

\indent Let $ n, r \in \mathbb{N} $ and let $ d_{1}, \ldots, d_{r} \in \mathbb{N} $ such that $ d_{1} \leq \ldots \leq d_{r} $. 

\begin{defn}
A \emph{polynomial vector} $ f= (f_1, \ldots, f_r) $ is a vector of $ r $ forms in $ n+1 $ variables over $ \mathbb{F} $ of degrees $ d_{1}, \ldots, d_{r} $, i.e. $ f_{j} \in \mathbb{F}[x_{0}, \ldots, x_{n}]_{d_{j}} $ for $ j \in \{1, \ldots, r\} $.  
\end{defn}

\begin{defn}\label{def:simWardec}
A \emph{simultaneous Waring decomposition} over $ \mathbb{F} $ of $ f= (f_1, \ldots, f_r) $ is given by $ \ell_{1}, \ldots, \ell_{k} \in \mathbb{F}[x_{0}, \ldots, x_{n}]_{1} $ and $ (\lambda_{1}^{j}, \ldots, \lambda_{k}^{j}) \in \mathbb{F}^{k} - \{\underline{0}\} $ s.t.
\begin{equation}\label{eq:simWardec}
f_{j} = \lambda_{1}^{j}\ell_{1}^{d_{j}} + \ldots + \lambda_{k}^{j}\ell_{k}^{d_{j}} 
\end{equation}
for all $ j \in \{1, \ldots, r\} $, or in vector notation
\begin{equation}\label{eq:vecWardec}
f = \sum_{i=1}^{k}\left(\lambda_{i}^{1}\ell_{i}^{d_{1}}, \ldots, \lambda_{i}^{r}\ell_{i}^{d_{r}}\right).
\end{equation}
\end{defn}

\begin{rem}
The adjective \emph{simultaneous} in Definition \ref{def:simWardec} is justified by the fact that $ f_{1}, \ldots, f_{r} $ are decomposed in the sense of Definition \ref{def:Wardec} by means of the same linear forms $ \ell_{1}, \ldots, \ell_{k} $. For simplicity, throughout the paper, we speak about Waring decomposition also for a polynomial vector.
\end{rem}

\begin{rem}
The notions introduced in section \ref{ss:dr} extend in a natural way to the case $ r > 1 $. In particular, referring to Definition \ref{def:rank}, Definition \ref{def:subgen} and Definition \ref{def:ident}, we say that:
\begin{itemize}
\item the \emph{Waring} \emph{rank} over $ \mathbb{F} $ of a polynomial vector $ f = (f_{1}, \ldots, f_{r}) $  is the minimal $ k $ appearing in (\ref{eq:vecWardec}); 
\item polynomial vectors over $ \mathbb{C} $ have \emph{sub-generic rank} $ k $ if 
\begin{equation}\label{eq:subgen}
k < \left\lceil{\frac{1}{n+r}} \sum_{j=1}^{r}{n+d_{j} \choose d_{j}}\right\rceil; 
\end{equation}
\item a rank-$ k $ polynomial vector $ f = (f_{1}, \ldots, f_{r}) $ is \emph{Waring identifiable} over $ \mathbb{F} $ if the expression (\ref{eq:vecWardec}) is unique up to reordering and rescaling.
\end{itemize}
\end{rem}

Unlike what occurs for the case $ r = 1 $ (section \ref{sss:ci}), in the simultaneous Waring setting, up to now, there are results only on the complex side and it is still an open problem to determine the complete list of identifiable cases of generic rank over $ \mathbb{C} $. To the best of our knowledge, the discovered ones are essentially $ 5 $ and they are listed in Table $ 4 $. 

\begin{table}[h]
\begin{center}
\begin{tabular} {c |c |c |c |c }
{\bf n} & {\bf r} &  ${\bf d_{1}, \ldots, d_{r}}$ & {\bf k}  & {\bf Ref.} \\
$ 1 $ & $ \forall $ & $ d_1+1 \geq k  $  & $ \left \lceil \frac{1}{1+r} \sum_{j=1}^{r}{1+d_{j} \choose d_{j}} \right \rceil $ & \cite{CR} \\
$ \forall $ & $2$ & $2,2$ & $ n+1 $ & \cite{We} \\
$ 2 $ & $2$ & $2,3$ & $ 4 $ & \cite{Ro} \\
$ 2 $ & $3$ & $3,3,4$ & $ 7 $ & \cite{AGMO} \\
$ 2 $ & $4$ & $2,2,2,2$ & $ 4 $ & \cite{Ve} \\
\end{tabular}\vspace{0.2cm}\caption{Simult. Waring identifiable cases of generic rank over $ \mathbb{C} $}
\end{center}
\end{table}

Concerning the second line of Table $ 4 $, which we mentioned in the Introduction, we stress that, while any quadric of generic rank $ n+1 $ has infinitely many minimal Waring decompositions over $ \mathbb{C} $, it becomes identifiable by requiring the simultaneous decomposition with another quadric. 

Similarly, the third line of Table $ 4 $ reveals that any cubic of generic rank $ 4 $, which is not Waring identifiable over $ \mathbb{C} $, admits a unique minimal presentation by adding a quadric. This case has a geometric interpretation, close to the previous one, by changing an orthogonal basis with a \emph{Parseval frame}, in the sense of \cite{CMS}. Indeed, assume that $ f = (f_{1},f_{2}) $ is a polynomial vector with $ f_{1} = x_{0}^2+x_{1}^2+x_{2}^{2} $ and $ f_{2} $ general in $ \mathbb{C}[x_{0},x_{1}, x_{2}]_{3} $. The unique Waring presentation of $ f $ consists of $ \ell_{1}, \ldots,\ell_{4} \in \mathbb{C}[x_{0},x_{1}, x_{2}]_{1} $ and $ (\lambda_{1}, \ldots, \lambda_{4}) \in \mathbb{C}^{4} - \{\underline{0}\} $ s.t.
$$
(f_{1},f_{2}) = \sum_{i=1}^{4}\left(\ell_{i}^{2}, \lambda_{i}\ell_{i}^{3}\right).
$$
The $ 3 \times 4 $ matrix $ L $, whose columns contain the coefficient of $ \ell_{1}, \ldots, \ell_{4} $ in the standard basis $\{x_{0}, x_{1},x_{2}\}$, satisfies $ L \overline{L}^{t}= I_{3} $, i.e. the rows of $ L $ form an orthonormal set of $  \mathbb{C}^{4} $. The columns of $ L $ yield a Parseval frame of $ \mathbb{C}^{3} $.

Therefore, introducing the simultaneous Waring decomposition over $ \mathbb{C} $ of more forms with the same (number of) summands, provides another useful technique to recover classical Waring identifiability, when it fails over $ \mathbb{C} $.

\section{Recent results}

\subsection{The generic rank case through Numerical Algebraic Geometry}\label{ss:NAG} $ \quad $ \\

By combining the points of view described in section \ref{sss:ri} and section \ref{sss:sci} we get the following:

\begin{thm}\label{thm:A}
Let $ {\mathbb{P}_{\mathbb{R}}^{35}} $ be the projective space over $ \mathbb{R} $ defined by real polynomial vectors with $ n = 2 $, $ r= 4 $ and $ d_1 = 2, d_2 = d_3 = d_4  = 3 $. \\
There exists a non-trivial Euclidean open set $ U \subset {\mathbb{P}_{\mathbb{R}}^{35}} $ such that any $ f  \in U $ has rank $ 6 $ over $ \mathbb{C} $ (and over $ \mathbb{R} $ too) and is identifiable over $ \mathbb{R} $ but not over $ \mathbb{C} $.
\end{thm} 

Theorem \ref{thm:A} arises from a computational approach via the software Bertini for \emph{Numerical Algebraic Geometry}. In this sense, in section \ref{sss:cap} we describe the computer-aided procedure that we implemented to get our result and in section \ref{sss:proof} we outline its proof.

\subsubsection{Computational technique}\label{sss:cap} $ \quad $ \\

We consider the polynomial system, equivalent to (\ref{eq:vecWardec}), with the assumptions $ n=2, r=4, d_{1}=2, d_{2}=d_{3}=d_{4}=3, k=6 $:

\begin{equation}\label{eq:polsys}
\left\{
\begin{array}{l}
f_{1} - \lambda_{1}^{1}\ell_{1}^{2}- \ldots -  \lambda_{6}^{1}\ell_{6}^{2} = 0 \\
f_{2} - \lambda_{1}^{2}\ell_{1}^{3}- \ldots -  \lambda_{6}^{2}\ell_{6}^{3} = 0 \\
f_{3} - \lambda_{1}^{3}\ell_{1}^{3}- \ldots -  \lambda_{6}^{3}\ell_{6}^{3} = 0 \\
f_{4} - \lambda_{1}^{4}\ell_{1}^{3}- \ldots -  \lambda_{6}^{4}\ell_{6}^{3} = 0 \\
\end{array}
\right.
\end{equation}

In (\ref{eq:polsys}), $ f_{j} \in \mathbb{R}[x_{0}, x_{1}, x_{2}]_{d_{j}} $ is a fixed general form, while $ \ell_{i} = x_{0}+ l_{1}^{i}x_{1}+ l_{2}^{i}x_{2} \in \C[x_{0}, x_{1}, x_{2}]_{1}$ and $ \lambda_{i}^{j} \in \C $ are unknown. 

By applying the identity principle for polynomials, the $j$-th equation of (\ref{eq:polsys}) splits in $ {{d_{j}+2} \choose {2}} $ conditions, so that we get a square non linear system of order $ 36$, which we denote by $ F_{\left(f_{1}, f_{2}, f_{3}, f_{4}\right)}\left(\left[l_{1}^{1}, l_{2}^{1}, \lambda_{1}^{1}, \lambda_{1}^{2},\lambda_{1}^{3},\lambda_{1}^{4}\right],\ldots, \left[l_{1}^{6}, l_{2}^{6}, \lambda_{6}^{1}, \lambda_{6}^{2},\lambda_{6}^{3},\lambda_{6}^{4}\right]\right) $, which is equivalent to (\ref{eq:polsys}). 

We aim to compute the number of solutions of $ F_{\left(f_{1}, f_{2}, f_{3}, f_{4}\right)} $ in $ \mathbb{R} $, that is the number of Waring decompositions over $ \mathbb{R} $ with $ 6 $ summands of the polynomial vector $ f = (f_{1}, f_{2}, f_{3}, f_{4}) $. 

In practice, to work with general $ f_{j} $'s, we assign random values $ \overline{l}_{h}^{i} $, $ \overline{\lambda}_{i}^{j} \in \mathbb{R} $ to $ l_{h}^{i} $, $ \lambda_{i}^{j} $, getting a real vector $ \left(\left[\overline{l}_{1}^{1}, \overline{l}_{2}^{1}, \overline{\lambda}_{1}^{1},\overline{\lambda}_{1}^{2},\overline{\lambda}_{1}^{3}, \overline{\lambda}_{1}^{4}\right], \ldots, \left[\overline{l}_{1}^{6}, \overline{l}_{2}^{6}, \overline{\lambda}_{6}^{1}, \overline{\lambda}_{6}^{2},\overline{\lambda}_{6}^{3},\overline{\lambda}_{6}^{4}\right]\right) \in \mathbb{R}^{36}$, called a \emph{start-point}. By means of $ F_{\left(f_{1}, f_{2}, f_{3}, f_{4}\right)} $, we compute the corresponding $ \overline{f}_{1}, \overline{f}_{2},\overline{f}_{3}, \overline{f}_{4} $, whose coefficients are called \emph{start-parameters}. The input of our procedure consists of start-point and start-parameters. In particular, the start-point is, by construction, a real solution of $ F_{\left(\overline{f}_{1},\overline{f}_{2},\overline{f}_{3}, \overline{f}_{4}\right)} $, i.e. it provides a Waring decomposition over $ \mathbb{R} $  with $ 6 $ summands of $ \overline{f} = \left(\overline{f}_{1}, \overline{f}_{2},\overline{f}_{3}, \overline{f}_{4}\right) $. 

Afterwards, we start a \emph{triangle-loop}, i.e. a loop divided into $ 3 $ steps. 

First, in $ F_{\left(\overline{f}_{1},\overline{f}_{2},\overline{f}_{3}, \overline{f}_{4}\right)} $ we replace the start-parameters with random complex entries, called \emph{final-parameters} and we get a new square polynomial system, $ F_{1} $, of order $36 $. We construct a segment homotopy
$$ H_{0} : \C^{36} \times [0,1] \to \C^{36} $$
between $ F_{\left(\overline{f}_{1}, \overline{f}_{2},\overline{f}_{3}, \overline{f}_{4}\right)} $ and $ F_{1} $: it provides a \emph{path} connecting the start-point to a solution of $ F_{1} $, called \emph{end-point}. 

The final-parameters and the end-point become, respectively, the start-parameters and the start-point for the second step. Starting from $ F_{1} $ and proceeding as in the first step, we obtain a polynomial system $ F_{2} $ and through a segment homotopy
$$ H_{1} : \C^{36} \times [0,1] \to \C^{36} $$
between $ F_{1} $ and $ F_{2} $, the start-point is sent to a solution of $ F_{2} $. 

Finally, in $ F_{2} $ the start-parameters are replaced with the entries of $ \overline{f}_{1}, \ldots, \overline{f}_{4} $, getting again $ F_{\left(\overline{f}_{1}, \overline{f}_{2},\overline{f}_{3}, \overline{f}_{4}\right)} $. A segment homotopy
$$ H_{2} : \C^{36} \times [0,1] \to \C^{36} $$
from $ F_{2} $ to $ F_{\left(\overline{f}_{1}, \overline{f}_{2},\overline{f}_{3}, \overline{f}_{4}\right)} $ sends the start-point of the third step to a solution of $ F_{\left(\overline{f}_{1}, \overline{f}_{2},\overline{f}_{3}, \overline{f}_{4}\right)} $, which is the output of our algorithm.

At the end of the triangle-loop, we check if the output differs from the initial start-point $ \left(\left[\overline{l}_{1}^{1}, \overline{l}_{2}^{1}, \overline{\lambda}_{1}^{1},\overline{\lambda}_{1}^{2},\overline{\lambda}_{1}^{3}, \overline{\lambda}_{1}^{4}\right], \ldots, \left[\overline{l}_{1}^{6}, \overline{l}_{2}^{6}, \overline{\lambda}_{6}^{1}, \overline{\lambda}_{6}^{2},\overline{\lambda}_{6}^{3},\overline{\lambda}_{6}^{4}\right]\right)$. If this is not the case, then we restart the procedure. Otherwise: if the output is complex but not real, then this procedure suggests that the polynomial vector $ \overline{f} = \left(\overline{f}_{1}, \overline{f}_{2},\overline{f}_{3}, \overline{f}_{4} \right) $ under investigation is real identifiable, while if the output is real, then identifiability over $ \mathbb{R} $ fails. 

\begin{rem}
This computational technique provides always an answer, since, as we will see in the forthcoming paper \cite{ABC1}, the general polynomial vector $ f = (f_{1},f_{2},f_{3},f_{4}) $ with $ f_{1} \in \mathbb{C}[x_{0},x_{1},x_{2}]_{2} $ and $ f_{j} \in \mathbb{C}[x_{0},x_{1},x_{2}]_{3} $, $ j \in \{2,3,4\} $, has $ 2 $ Waring decompositions over $ \mathbb{C} $ with rank $ 6 $. This fact had previously been checked by means of the computational analysis developed for the \emph{perfect cases} in \cite{AGMO}, from which section \ref{sss:cap} is inspired.
\end{rem}

\subsubsection{Proof of Theorem \ref{thm:A}}\label{sss:proof} $ \quad $ \\

We apply the method described in section \ref{sss:cap} to the real polynomial vector $ \overline{f} = (\overline{f}_{1}, \overline{f}_{2},\overline{f}_{3}, \overline{f}_{4}) $ arising from (\ref{eq:polsys}) with start-point

{\small{$$  \left[\overline{l}_{1}^{1}, \overline{l}_{2}^{1}, \overline{\lambda}_{1}^{1},\overline{\lambda}_{1}^{2},\overline{\lambda}_{1}^{3}, \overline{\lambda}_{1}^{4}\right] = \left[-0.89, -0.38, 0.38, 0.87, -0.50, 0.44 \right] \, $$
$$  \,\,\,\,\left[\overline{l}_{1}^{2}, \overline{l}_{2}^{2}, \overline{\lambda}_{2}^{1},\overline{\lambda}_{2}^{2},\overline{\lambda}_{2}^{3}, \overline{\lambda}_{2}^{4}\right]= \left[0.71, -0.46, -0.54, -0.62, 0.22,-0.37 \right]$$
$$  \,\,\,\, \left[\overline{l}_{1}^{3}, \overline{l}_{2}^{3}, \overline{\lambda}_{3}^{1},\overline{\lambda}_{3}^{2},\overline{\lambda}_{3}^{3}, \overline{\lambda}_{3}^{4}\right] = \left[0.88, -0.50, -0.92, 0.86, -0.74, -0.74 \right] $$
$$  \left[\overline{l}_{1}^{4}, \overline{l}_{2}^{4}, \overline{\lambda}_{4}^{1},\overline{\lambda}_{4}^{2},\overline{\lambda}_{4}^{3}, \overline{\lambda}_{4}^{4}\right] = \left[-0.50, 0.72, 0.73, 0.93, -0.30, -0.45 \right] \, $$
$$  \left[\overline{l}_{1}^{5}, \overline{l}_{2}^{5}, \overline{\lambda}_{5}^{1},\overline{\lambda}_{5}^{2},\overline{\lambda}_{5}^{3}, \overline{\lambda}_{5}^{4}\right]= \left[0.39, 0.32, 0.16, 0.74, -0.78, 0.64 \right] \,\,\,\,\,\,\,\,\,\, $$
$$  \,\,\,\,\,\, \left[\overline{l}_{1}^{6}, \overline{l}_{2}^{6}, \overline{\lambda}_{6}^{1},\overline{\lambda}_{6}^{2},\overline{\lambda}_{6}^{3}, \overline{\lambda}_{6}^{4}\right] = \left[-0.52, -0.55, -0.99, -0.57, 0.68, 0.98 \right]. $$ }}
For simplicity, the start-point has been divided into $ 6 $ blocks, corresponding to the $ 6 $ summands in the Waring decomposition of $ \overline{f} $. In particular, in the each block, the first $ 2 $ entries are the coefficients of the linear form, while the other $ 4 $ are the scalars appearing in front of the power of this linear form in the decomposition.

The solutions of $ F_{(\overline{f}_{1}, \overline{f}_{2},\overline{f}_{3}, \overline{f}_{4})} $ are the start-point and the complex but not real point whose $ 6 $ blocks are the following:

{\small{
$$  \left[\underline{l}_{1}^{1}, \underline{l}_{2}^{1}, \underline{\lambda}_{1}^{1},\underline{\lambda}_{1}^{2},\underline{\lambda}_{1}^{3}, \underline{\lambda}_{1}^{4}\right] = \left[-6.783965763899463\cdot10^{-1}-i\,3.078910418301080\cdot10^{-1}, \right. $$
$$ \quad\quad\quad\quad\quad\quad\quad\quad\quad\quad -3.934665356579067\cdot10^{-1}-i\,6.002384989512501\cdot10^{-2},$$
$$\quad\quad\quad\quad\quad\quad\quad\quad\quad\quad 1.635655698401628\cdot10^{-1}-i\, 3.203959975350376\cdot10^{-1},$$
$$\quad\quad\quad\quad\quad\quad\quad\quad\quad\quad 4.464358529188056\cdot10^{-1}-i \, 6.267867847785145\cdot10^{-1},$$
$$\quad\quad\quad\quad\quad\quad\quad\quad\quad\quad -2.490818000203756\cdot10^{-1}+i\,3.764670668431170\cdot10^{-1},$$
$$\quad\quad\quad\quad\quad\quad\quad\quad\quad\quad \left. 2.681564395007827\cdot10^{-1}-i\,1.891499636461206\cdot10^{-1}\right] $$
$$  \,\,\,\left[\underline{l}_{1}^{2}, \underline{l}_{2}^{2}, \underline{\lambda}_{2}^{1},\underline{\lambda}_{2}^{2},\underline{\lambda}_{2}^{3}, \underline{\lambda}_{2}^{4}\right]= \left[-4.772309397135773\cdot10^{-1}-i\,1.199357995730621\cdot10^{-15}, \right.$$
$$ \quad\quad\quad\quad\quad\quad\quad\quad\quad\quad 6.992780172912881\cdot10^{-1}+i\,2.400152559339785\cdot10^{-17},$$
$$\quad\quad\quad\quad\quad\quad\quad\quad\quad\quad 7.835648558121242\cdot10^{-1}-i\, 4.895172808772053\cdot10^{-17},$$
$$\quad\quad\quad\quad\quad\quad\quad\quad\quad\quad 9.966810767997231\cdot10^{-1}+i \, 2.151057110211241\cdot10^{-16},$$
$$\quad\quad\quad\quad\quad\quad\quad\quad\quad\quad -3.200804961183177\cdot10^{-1}-i\,7.470017443153565\cdot10^{-16},$$
$$\quad\quad\quad\quad\quad\quad\quad\quad\quad\quad \left. -4.673109022088685\cdot10^{-1}+i\,1.011357339021635\cdot10^{-15} \right] $$
$$  \,\left[\underline{l}_{1}^{3}, \underline{l}_{2}^{3}, \underline{\lambda}_{3}^{1},\underline{\lambda}_{3}^{2},\underline{\lambda}_{3}^{3}, \underline{\lambda}_{3}^{4}\right] = \left[3.684954187786456\cdot10^{-1}-i\,4.971134723481818\cdot10^{-16}, \right.$$
$$ \quad\quad\quad\quad\quad\quad\quad\quad\quad\quad 3.097626827511096\cdot10^{-1}-i\,1.013037852388987\cdot10^{-15},$$
$$\quad\quad\quad\quad\quad\quad\quad\, 1.203175963086619+i\, 1.411261922211124\cdot10^{-15},$$
$$\quad\quad\quad\quad\quad\quad\quad\quad\quad\quad 7.178706316282033\cdot10^{-1}+i \, 2.830757004669560\cdot10^{-15},$$
$$\quad\quad\quad\quad\quad\quad\quad\quad\quad\quad -8.000652140416307\cdot10^{-1}-i\,1.314324083595553\cdot10^{-15},$$
$$\quad\quad\quad\quad\quad\quad\quad\, \left. 6.688294919987063+i\,6.432029388270255\cdot10^{-17}\right] $$
$$ \,\,\,\,\,\, \left[\underline{l}_{1}^{4}, \underline{l}_{2}^{4}, \underline{\lambda}_{4}^{1},\underline{\lambda}_{4}^{2},\underline{\lambda}_{4}^{3}, \underline{\lambda}_{4}^{4}\right] = \left[-5.642110773943350\cdot10^{-1}+i\,7.105047886840632\cdot10^{-16}, \right.$$
$$ \quad\quad\quad\quad\quad\quad\quad\quad\quad\quad -5.613848544445357\cdot10^{-1}-i\,7.586704701967317\cdot10^{-16},$$
$$\quad\quad\quad\quad\quad\quad\quad\, -1.049596368398710-i\, 5.833007687972014\cdot10^{-15},$$
$$\quad\quad\quad\quad\quad\quad\quad\quad\quad\quad -8.493748090457428\cdot10^{-1}-i \,1.022155992659590\cdot10^{-14},$$
$$\quad\quad\quad\quad\quad\quad\quad\quad\quad\quad 8.143985201225270\cdot10^{-1}+i\,4.665972469508617\cdot10^{-15},$$
$$\quad\quad\quad\quad\quad\quad\quad\quad\quad\quad \left. 7.945024318077690\cdot10^{-1}-i\,7.342420399979177\cdot10^{-15}\right] $$
$$  \left[\underline{l}_{1}^{5}, \underline{l}_{2}^{5}, \underline{\lambda}_{5}^{1},\underline{\lambda}_{5}^{2},\underline{\lambda}_{5}^{3}, \underline{\lambda}_{5}^{4}\right]=\left[-6.783965763898687\cdot10^{-1}+i\,3.078910418302253\cdot10^{-1}, \right.$$
$$ \quad\quad\quad\quad\quad\quad\quad\quad\quad\quad -3.934665356579083\cdot10^{-1}+i\,6.002384989511068\cdot10^{-2},$$
$$\quad\quad\quad\quad\quad\quad\quad\quad\quad\quad 1.635655698401061\cdot10^{-1}+i\, 3.203959975350371\cdot10^{-1},$$
$$\quad\quad\quad\quad\quad\quad\quad\quad\quad\quad 4.464358529186822\cdot10^{-1}+i\,6.267867847785181\cdot10^{-1},$$
$$\quad\quad\quad\quad\quad\quad\quad\quad\quad\quad -2.490818000203047\cdot10^{-1}-i\,3.764670668431184\cdot10^{-1},$$
$$\quad\quad\quad\quad\quad\quad\quad\quad\quad\quad \left. 2.681564395007271\cdot10^{-1}+i\,1.891499636461248\cdot10^{-1} \right] $$
$$  \left[\underline{l}_{1}^{6}, \underline{l}_{2}^{6}, \underline{\lambda}_{6}^{1},\underline{\lambda}_{6}^{2},\underline{\lambda}_{6}^{3}, \underline{\lambda}_{6}^{4}\right] = \left[8.537942728260280\cdot10^{-1}+i\,1.366691048526603\cdot10^{-15}, \right.$$
$$ \quad\quad\quad\quad\quad\quad\quad\quad\quad\quad -4.980478154684624\cdot10^{-1}+i\,9.271283827466670\cdot10^{-17},$$
$$\quad\quad\quad\quad\quad\quad\quad\quad\quad\quad -1.361417223402345+i\, 5.131989668297199\cdot10^{-15},\,\,\,\,\,\,\,\,\,\,\,\,\,\,\,\,\,$$
$$\quad\quad\quad\quad\quad\quad\quad\quad\quad\quad 4.519513947803280\cdot10^{-1}+i\,3.420576539028830\cdot10^{-15},$$
$$\quad\quad\quad\quad\quad\quad\quad\quad\quad\quad -6.160892099218980\cdot10^{-1}-i\,1.042175785774535\cdot10^{-15},$$
$$\quad\quad\quad\quad\quad\quad\quad\quad\quad\quad \left. -1.032333900599116 +i\,1.950994704333353\cdot10^{-15}\right].\,\,\,\,\,\,\,\,\,\,\,\,\, $$
}}

\noindent We notice that $ \left(\left[\underline{l}_{1}^{1}, \underline{l}_{2}^{1}, \underline{\lambda}_{1}^{1},\underline{\lambda}_{1}^{2},\underline{\lambda}_{1}^{3}, \underline{\lambda}_{1}^{4}\right], \ldots, \left[\underline{l}_{1}^{6}, \underline{l}_{2}^{6}, \underline{\lambda}_{6}^{1}, \underline{\lambda}_{6}^{2},\underline{\lambda}_{6}^{3},\underline{\lambda}_{6}^{4}\right]\right)$ is \emph{self-conjugate}, i.e. it is a fixed point of the conjugate operator. Namely, by considering rounding errors, the summands corresponding to the first and the fifth block are one the conjugate of the other and the remaining $ 4 $ are real.

If we move $\overline{f}=(\overline{f}_{1}, \overline{f}_{2},\overline{f}_{3}, \overline{f}_{4})$ in a small Euclidean open subset over $ \mathbb{R} $, only one decomposition remains real, because the property of being real is open in the set of Waring decompositions. Therefore there exists a non-trivial Euclidean open subset $ U $ in the space of real polynomial vectors under investigation whose elements have one real decomposition plus a complex not real one, as claimed.

\begin{rem}
A result similar to the one stated in Theorem \ref{thm:A} may be true for the cases of general real polynomial vectors with finitely many Waring decompositions over $ \mathbb{C} $. The interesting fact is that this phenomenon occurs also when the number of decompositions is even. 
\end{rem}

\subsection{The sub-generic rank case through the Hessian criterion}\label{sec:hes} $ \quad $ \\

A sufficient condition for Waring identifiability over $ \mathbb{C} $ of general polynomial vectors, whose rank satisfies (\ref{eq:subgen}), is given by the \emph{Hessian criterion}, which is based on the following generalization of Lemma 5.1 of \cite{COV1} concerning the case $ r = 1 $.

\begin{lem}[Sufficient condition for identifiability of general elements]\label{lemma:hess}$\quad$\\
Assume that the variety $ X^{n}_{d_{1},\ldots,d_{r}} = \mathbb{P}(\bigoplus_{j=1}^{r}\mathcal{O}_{\mathbb{P}_{\mathbb{C}}^n}(d_{j})) $ of rank-$ 1 $ complex polynomial vectors whose $ r $ components depend on $ n+1 $ variables and have degrees $ d_{1} \leq \ldots \leq d_{r} $ is not $ k $-defective, with $ k $ as in (\ref{eq:subgen}). Let $ g_{1}, \ldots, g_{k} \in X^{n}_{d_{1},\ldots,d_{r}} $ be general polynomial vectors and $ f \in \langle g_{1}, \ldots, g_{k}\rangle $ general. Let $ T = \langle T_{g_{1}}X^{n}_{d_{1},..,d_{r}},\ldots, T_{g_{k}}X^{n}_{d_{1},..,d_{r}}\rangle $ and let $ \mathcal{C}_{k} = \{g \in X^{n}_{d_{1},\ldots,d_{r}} \, | \, T_{g}X^{n}_{d_{1},..,d_{r}} \subset T \} $ be the $ k $-tangential contact locus.
If
\begin{itemize}
\item[$ \bullet $] $ T $ has the expected dimension $ k(\dim X^{n}_{d_{1},..,d_{r}}+1) = k(n+r) $
\item[$ \bullet $] $ \mathcal{C}_{k} $ is $ 0 $-dimensional at each $ g_{i} $
\end{itemize}
then $ f $ is Waring identifiable over $ \mathbb{C} $ with rank $ k $. 
\end{lem} 

In the following section \ref{sss:hp} we describe the algorithm, based on Lemma \ref{lemma:hess}, that we implemented via Macaulay2 to detect identifiable cases of sub-generic rank. With this procedure we are able to produce many identifiable cases, which are listed in section \ref{sss:res}.

\subsubsection{Computational technique}\label{sss:hp} $ \quad $ \\

After establishing the initial data $ n,r,d_{1}, \ldots,d_{r}$ and assigning to $ k $ the first value that satisfies (\ref{eq:subgen}), we define a parametrization of $ X^{n}_{d_{1},\ldots,d_{r}} \subset \mathbb{P}^{N-1} $, where $ N = \sum_{j=1}^{r}{n+d_{j} \choose d_{j}} $.

Then, we choose $ k $ random points $ g_{1}, \ldots, g_{k} \in X^{n}_{d_{1},\ldots,d_{r}} $, which guarantees that they are general, and we construct the corresponding $ k(n+r)\times N $ matrix of constants $ j_1 $ associated to $ T $. If the rank of $ j_{1} $ equals $ k(n+r) $, then $ T $ has the expected dimension (and $ X^{n}_{d_{1},\ldots,d_{r}} $ is not $ k $-defective because of Terracini's Lemma), otherwise we decrease $ k $ by $ 1 $ and we go back to the choice of points.

By computing the kernel of $ j_{1} $ we obtain the $ N-k(n+r) $ Cartesian equations of $ T $ and by partial derivation of them we get the $ (n+r)\times (N-k(n+r)) $ Cartesian equations of $ \mathcal{C}_{k}. $ 

In order to determine the dimension of $ \mathcal{C}_k $ at $ g_i $, $ i \in \{1, \ldots, k \} $, we compute the codimension of $ T_{g_i}\mathcal{C}_{k} $. It equals the rank of the Jacobian, evaluated at $ g_{i} $, of the Cartesian equations of $ \mathcal{C}_{k} $, i.e. of a matrix $ H $, $ (n+r)\times((n+r)\times (N-k(n+r))) $, with constant entries. This procedure is an equivalent way to construct the \emph{Hessian} matrix of the equations defining $ T $, which justifies the expression Hessian criterion. If $ \rank H = \dim X^{n}_{d_{1},..,d_{r}} = n+r-1 $, then $ \mathcal{C}_{k} $ has dimension $ 0 $ at $ g_{i} $ and, by Lemma \ref{lemma:hess}, the general complex $ f = (f_{1}, \ldots, f_{r}) \in \langle g_{1}, \ldots, g_{k}\rangle $ of rank $ k $ is Waring identifiable over $ \mathbb{C} $. Otherwise, other techniques are needed to confirm or deny the identifiability for this particular value of the rank $ k $. Thus, we decrease $ k $ by $ 1 $ and we restart the procedure from the choice of the $ k $ points.

As output of this algorithm, we get the highest value of the sub-generic rank for which the Hessian criterion implies identifiability over $ \mathbb{C} $ (note that, if identifiability holds for certain $ \tilde{k} $, then the same is true for any $ k \leq \tilde{k} $, \cite{CO}).

\subsubsection{Table of results}\label{sss:res} $ \quad $ \\

In the following table, we list the results obtained via the computational technique described in section \ref{sss:hp}. In the heading of the third column, we denote by $ g $ the expected generic rank, given by the formula $ g = \left \lceil{\frac{1}{n+r}}\sum_{j=1}^{r}{n+d_{j} \choose d_{j}} \right \rceil $. In particular, we adopt the notation $ g^{\star} $ for the defective cases, with the corresponding bibliographic reference, and $ g^{\diamond} $ for the cases for whom the defectivity problem is still open.
$ \quad $
\begin{longtable}{|c|c|l|c|c|c|}
\hline \multicolumn{1}{|c|}{${\bf r}$} & \multicolumn{1}{|c|}{${\bf n}$}& \multicolumn{1}{|l|}{ ${\bf (d_{1},\ldots,d_{r})}$}& \multicolumn{1}{|c|}{${\bf g}$}& \multicolumn{1}{|c|}{${\bf k}$}\\
\hline
\endhead
\hline
\endfoot
$2$ &$2$ &$(3,3)$ &5$^{\star}$\cite{CC} & $ \leq4 $ \\
$2$ &$2$ &$(4,4)$ & 8 & $ \leq7 $ \\
$2$ &$2$ &$(5,5)$ & 11 & $ \leq10 $ \\
$2$ &$2$ &$(6,6)$ & 14 &$ \leq 13 $  \\
$2$ &$2$ &$(7,7)$ & 18 & $ \leq 17 $  \\
$2$ &$2$ &$(8,8)$ & 23& $ \leq 22 $  \\
$2$ &$2$ &$(9,9)$ & 28 & $ \leq 27 $  \\
$2$ &$2$ &$(10,10)$ & 33 & $ \leq 32 $  \\
$2$ &$2$ &$(2,4)$ & 6&$ \leq4 $\\
$2$ &$2$ &$(2,5)$ & 7& $ \leq5 $\\
$2$ &$2$ &$(2,6)$ & 9$^{\diamond}$& $ \leq5 $\\
$2$ &$2$ &$(2,7)$ & 11$^{\diamond}$&$ \leq5 $\\
$2$ &$2$ &$(2,8)$ & 13$^{\diamond}$ &$ \leq5 $\\
$2$ &$2$ &$(2,9)$ & 16$^{\diamond}$& $ \leq5 $\\
$2$ &$2$ &$(2,10)$ & 18$^{\diamond}$& $ \leq5 $\\
$2$ &$2$ &$(3,4)$ & 7& $ \leq6 $\\
$2$ &$2$ &$(3,5)$ & 8& $ \leq7 $\\
$2$ &$2$ &$(3,6)$ & 10& $ \leq8 $\\
$2$ &$2$ &$(3,7)$ & 12&$ \leq9 $\\
$2$ &$2$ &$(3,8)$ &  14$^{\diamond}$& $ \leq8 $\\
$2$ &$2$ &$(3,9)$ &  17$^{\diamond}$& $ \leq9 $\\
$2$ &$2$ &$(3,10)$ &  19$^{\diamond}$& $ \leq9 $\\
$2$ &$2$ &$(4,5)$ &  9& $ \leq8 $\\
$2$ &$2$ &$(4,6)$ &  11& $ \leq10 $\\
$2$ &$2$ &$(4,7)$ &  13& $ \leq12 $\\
$2$ &$2$ &$(4,8)$ &  15& $ \leq14 $\\
$2$ &$2$ &$(4,9)$ &  18$^{\diamond}$& $ \leq14 $\\
$2$ &$2$ &$(4,10)$ &  21$^{\diamond}$& $ \leq14 $\\
$2$ &$2$ &$(6,7)$ &  16& $ \leq15 $\\
$2$ &$2$ &$(6,8)$ &  19& $ \leq18 $\\
$2$ &$2$ &$(6,9)$ &  21& $ \leq20 $\\
$2$ &$2$ &$(6,12)$ &  30$^{\diamond}$& $ \leq27 $\\
$3$ &$2$ &$(2,2,2)$ &4& $ \leq3 $ \\
$3$ &$2$ &$(3,3,3)$ &6& $ \leq5 $ \\
$3$ &$2$ &$(4,4,4)$ & 9 & $ \leq8 $ \\
$3$ &$2$ &$(5,5,5)$ & 13 & $ \leq12 $ \\
$3$ &$2$ &$(6,6,6)$ & 17 &$ \leq 16 $  \\
$3$ &$2$ &$(7,7,7)$ & 22 & $ \leq 21 $  \\
$3$ &$2$ &$(8,8,8)$ & 27& $ \leq 26 $  \\
$3$ &$2$ &$(9,9,9)$ & 33 & $ \leq 32 $  \\
$3$ &$2$ &$(10,10,10)$ & 40 & $ \leq 39 $  \\
$3$ &$2$ &$(2,3,3)$ &6& $ \leq4 $ \\
$3$ &$2$ &$(2,4,4)$ &8& $ \leq5 $ \\
$3$ &$2$ &$(3,4,4)$ & 8 & $ \leq7 $ \\
$3$ &$2$ &$(3,5,5)$ & 11 & $ \leq9 $ \\
$3$ &$2$ &$(2,2,3)$ & 5 &$ \leq 3 $  \\
$3$ &$2$ &$(2,2,4)$ &6 & $ \leq 4 $  \\
$3$ &$2$ &$(2,2,5)$ & 7& $ \leq 5 $  \\
$3$ &$2$ &$(2,2,6)$ & 8$^{\diamond}$ & $ \leq 5 $  \\
$3$ &$2$ &$(3,3,4)$ & 7 & $ \leq 6 $  \\
$3$ &$2$ &$(3,3,5)$ & 9 & $ \leq 8 $  \\
$3$ &$2$ &$(3,3,6)$ & 10 & $ \leq 8 $  \\
$3$ &$2$ &$(3,3,7)$ & 12 & $ \leq 9 $  \\
$3$ &$2$ &$(4,4,5)$ & 11 & $ \leq 10 $  \\
$3$ &$2$ &$(4,4,6)$ & 12 & $ \leq 11 $  \\
$3$ &$2$ &$(4,4,7)$ & 14 & $ \leq 13 $  \\
$3$ &$2$ &$(4,4,8)$ & 15 & $ \leq 14 $  \\
$3$ &$2$ &$(4,4,9)$ & 17$^{\diamond}$ & $ \leq 14 $  \\
$3$ &$2$ &$(5,5,6)$ & 14 & $ \leq 13 $  \\
$3$ &$2$ &$(5,5,7)$ & 16 & $ \leq 15 $  \\
$3$ &$2$ &$(5,5,8)$ & 18 & $ \leq 17 $  \\
$3$ &$2$ &$(5,5,9)$ & 20& $ \leq 19 $  \\
$3$ &$2$ &$(5,5,10)$ & 22 & $ \leq 20 $  \\
$3$ &$2$ &$(6,6,7)$ & 19 & $ \leq 18 $  \\
$3$ &$2$ &$(6,6,8)$ & 21 & $ \leq 20 $  \\
$3$ &$2$ &$(6,6,9)$ & 23& $ \leq 22 $  \\
$3$ &$2$ &$(6,6,10)$ & 25 & $ \leq 24 $  \\
$3$ &$2$ &$(6,6,11)$ & 27 & $ \leq 26 $  \\
$3$ &$2$ &$(6,6,12)$ & 30$^{\diamond}$ & $ \leq 27 $  \\
$3$ &$2$ &$(2,3,4)$ &7& $ \leq5 $ \\
$3$ &$2$ &$(2,3,5)$ &8& $ \leq5 $ \\
$3$ &$2$ &$(2,3,6)$ &9$^{\diamond}$& $ \leq5 $ \\
$3$ &$2$ &$(2,4,5)$ &9$^{\diamond}$& $ \leq5 $ \\
$3$ &$2$ &$(2,4,6)$ &10$^{\diamond}$& $ \leq5 $ \\
$3$ &$2$ &$(3,4,5)$ & 10 & $ \leq8 $ \\
$3$ &$2$ &$(3,4,6)$ & 11 & $ \leq9 $ \\
$3$ &$2$ &$(3,4,7)$ & 13 & $ \leq9 $ \\
$3$ &$2$ &$(4,5,6)$ & 13 & $ \leq 12 $  \\
$3$ &$2$ &$(4,5,7)$ & 15 & $ \leq 14 $  \\
$3$ &$2$ &$(4,5,8)$ & 17 & $ \leq 14 $  \\
$3$ &$2$ &$(4,5,9)$ & 19 & $ \leq 14 $  \\
$4$ &$2$ &$(3,3,3,3)$ &7& $ \leq6 $ \\
$4$ &$2$ &$(4,4,4,4)$ & 10 & $ \leq9 $ \\
$4$ &$2$ &$(5,5,5,5)$ & 14 & $ \leq13 $ \\
$4$ &$2$ &$(2,3,3,3)$ &6& $ \leq5 $ \\
$4$ &$2$ &$(2,2,3,3)$ &6& $ \leq4 $ \\
$4$ &$2$ &$(2,2,4,4)$ &7$^{\diamond}$& $ \leq5 $ \\
$4$ &$2$ &$(2,2,2,3)$ &5$^{\diamond}$& $ \leq3 $ \\
$4$ &$2$ &$(2,2,2,4)$ &6& $ \leq4 $ \\
$4$ &$2$ &$(3,4,5,5)$ &12& $ \leq9 $ \\
$4$ &$2$ &$(3,4,4,5)$ &11& $ \leq9 $ \\
$4$ &$2$ &$(3,3,4,5)$ &10& $ \leq8 $ \\
$5$ &$2$ &$(2,2,2,2,2)$ &5& $ \leq3 $ \\
$5$ &$2$ &$(3,3,3,3,3)$ &8& $ \leq6 $ \\
$5$ &$2$ &$(4,4,4,4,4)$ & 11 & $ \leq10 $ \\
$5$ &$2$ &$(5,5,5,5,5)$ & 15 & $ \leq14 $ \\
$5$ &$2$ &$(2,3,3,3,3)$ &7& $ \leq5 $ \\
$5$ &$2$ &$(2,2,2,2,3)$ & 5$^{\diamond}$ & $ \leq3 $ \\
$5$ &$2$ &$(2,2,2,2,4)$ & 6 & $ \leq4 $ \\
$6$ &$2$ &$(2,\ldots,2)$ &5$^{\star}$\cite{AB} & $ \leq3 $ \\
$6$ &$2$ &$(3,\ldots,3)$ &8& $ \leq7 $ \\
$6$ &$2$ &$(4,\ldots,4)$ & 12 & $ \leq10 $ \\
$6$ &$2$ &$(5,\ldots,5)$ & 16 & $ \leq15 $ \\
$6$ &$2$ &$(2,3,\ldots,3)$ &7& $ \leq5 $ \\
$6$ &$2$ &$(2,\ldots,2,3)$ & 5$^{\diamond}$ & $ \leq3 $ \\
$6$ &$2$ &$(2,\ldots,2,4)$ & 6 & $ \leq4 $ \\
$7$ &$2$ &$(2,\ldots,2)$ &5$^{\star}$\cite{AB} & $ \leq3 $ \\
$7$ &$2$ &$(4,\ldots,4)$ & 12 & $ \leq11 $ \\
$7$ &$2$ &$(5,\ldots,5)$ & 17 & $ \leq16 $ \\
$7$ &$2$ &$(2,3,\ldots,3)$ &8& $ \leq5 $ \\
$7$ &$2$ &$(2,\ldots,2,3)$ & 6 & $ \leq3 $ \\
$8$ &$2$ &$(2,\ldots,2)$ &5$^{\star}$\cite{AB}& $ \leq3 $ \\
$8$ &$2$ &$(3,\ldots,3)$ &8& $ \leq7 $ \\
$8$ &$2$ &$(4,\ldots,4)$ & 12 & $ \leq11 $ \\
$8$ &$2$ &$(5,\ldots,5)$ & 17 & $ \leq16 $ \\
$8$ &$2$ &$(2,3,\ldots,3)$ &8& $ \leq5 $ \\
$9$ &$2$ &$(3,\ldots,3)$ &9& $ \leq7 $ \\
$9$ &$2$ &$(4,\ldots,4)$ & 13 & $ \leq11 $ \\
$9$ &$2$ &$(5,\ldots,5)$ & 18 & $ \leq16 $ \\
$10$ &$2$ &$(3,\ldots,3)$ &9$^{\star}$\cite{AB}& $ \leq7 $ \\
$10$ &$2$ &$(4,\ldots,4)$ & 13 & $ \leq12 $ \\
$10$ &$2$ &$(5,\ldots,5)$ & 18 & $ \leq17 $ \\
$11$ &$2$ &$(2,\ldots,2)$ &6& $ \leq3 $ \\
$11$ &$2$ &$(4,\ldots,4)$ & 13 & $ \leq12 $ \\
$11$ &$2$ &$(5,\ldots,5)$ & 18 & $ \leq17 $ \\
$12$ &$2$ &$(4,\ldots,4)$ & 13 & $ \leq12 $ \\
$12$ &$2$ &$(5,\ldots,5)$ & 18 & $ \leq17 $ \\
$13$ &$2$ &$(4,\ldots,4)$ & 13 & $ \leq12 $ \\
$13$ &$2$ &$(5,\ldots,5)$ & 19 & $ \leq17 $ \\
$14$ &$2$ &$(4,\ldots,4)$ & 14 & $ \leq12 $ \\
$14$ &$2$ &$(5,\ldots,5)$ & 19 & $ \leq17 $ \\
$15$ &$2$ &$(4,\dots,4)$ & 14$^{\star}$\cite{AB} & $ \leq12 $ \\
$15$ &$2$ &$(5,\ldots,5)$ & 19 & $ \leq18 $ \\
$16$ &$2$ &$(4,\ldots,4)$ & 14$^{\star}$\cite{AB} & $ \leq12 $ \\
$16$ &$2$ &$(5,\ldots,5)$ & 19 & $ \leq18 $ \\
$19$ &$2$ &$(5,\ldots,5)$ & 19 & $ \leq18 $ \\
$20$ &$2$ &$(5,\ldots,5)$ & 19 & $ \leq18 $ \\
$21$ &$2$ &$(5,\ldots,5)$ & 20$^{\star}$\cite{AB} & $ \leq18 $ \\
$3$ &$3$ &$(2,3,3)$ &9& $ \leq7 $ \\
$4$ &$3$ &$(2,3,3,3)$ & 10 & $ \leq9 $ \\
\end{longtable}

\vspace{0.5cm} \indent {\it
A\,c\,k\,n\,o\,w\,l\,e\,d\,g\,m\,e\,n\,t\,s.\;} This paper arises partially from the conference ``Waring decompositions and identifiability via Bertini and Macaulay2 softwares", given by the author during the fourteenth International conference on Effective Methods in Algebraic Geometry - MEGA 2017 held at the Castle of Nice University, campus Valrose, in Nice (France), June 12-16 2017. The author would like to thank all the members of the committees for this interesting and stimulating opportunity. This research was partially supported by the Italian GNSAGA-INDAM, by the Italian PRIN2015 - Geometry of Algebraic Varieties (B16J15002000005) and by the Italian ``Progetto strategico di ricerca di base Anno 2014 - prof. Giorgio Ottaviani" of the University of Florence.

\end{document}